\documentclass[12pt]{article}
\usepackage{amsmath,amsfonts,amssymb,amsthm}

\newtheorem*{theo}{Theorem}
\newtheorem{theor}{Theorem}
\newtheorem{prop}{Proposition}
\newtheorem*{coro}{Corollary}
\newtheorem*{lem}{Lemma}

\newtheorem{conj}{Conjecture}
\newtheorem*{konj}{Conjecture 1'}
\theoremstyle{remark}

\numberwithin{equation}{section}
\newcommand\la{\lambda}
\newcommand\al{\alpha}

\newcommand{\bm}[1]{{\mbox{\boldmath $#1$}}}
\newcommand{\bmp}{\bm{p}}
\newcommand{\bmq}{\bm{q}}

\title{A positivity conjecture for Jack polynomials}
\author{Michel Lassalle\\
\small Centre National de la Recherche Scientifique\\[-0.8ex]
\small Institut Gaspard-Monge, Universit\'e de Marne-la-Vall\'ee\\[-0.8ex]
\small 77454 Marne-la-Vall\'ee Cedex, France\\[-0.8ex]
\small \texttt{lassalle @ univ-mlv.fr}\\[-0.8ex]
\small \texttt{http://igm.univ-mlv.fr/{\textasciitilde}lassalle}}
\date{}

\begin{document}
\maketitle
\begin{abstract}
We present a positivity conjecture for the coefficients of the development of Jack polynomials in terms of power sums. This extends Stanley's ex-conjecture about normalized characters of the symmetric group. We prove this conjecture for partitions having a rectangular shape.
\end{abstract}

\section{Introduction}

A finite weakly decreasing sequence of positive integers $\la= (\la_1,...,\la_r)$ is called a partition with length $l(\la)=r$ and weight $|\la|=\sum_i \la_i$. Let $n$ be a fixed positive integer and $\mathfrak{S}_n$ the group of permutations of $n$ letters. The irreducible representations of $\mathfrak{S}_n$ are labelled by partitions with weight $n$. Let $\chi^\la$ denote the corresponding character.

By identification of the cycle decomposition of a permutation with a partition, any partition $\rho$ with weight $n$ defines a conjugacy class of $\mathfrak{S}_n$. Let $\mu$ be a partition with weight $k\leq n$, and $(\mu,1^{n-k})$ the partition obtained by adding $n-k$ unary parts. The normalized character $\widehat{\chi}^\lambda_{\mu,1^{n-k}}$ is defined by
\[\widehat{\chi}^\la_{\mu,1^{n-k}}=\frac{n!}{(n-k)!} \frac{
    \chi^\la_{\mu,1^{n-k}}}{\chi^\la_{1^n}}, \]
where $\chi^\la_{\rho}$ is the value of $\chi^\la$ on the conjugacy classs $\rho$, and $\chi^\la_{1^n}$ is the dimension of the irreducible representation.

Let $\bmp=(p_1,\ldots,p_m)$ and $\bmq=(q_1,\ldots,q_m)$, with
$q_1\geq \cdots \geq q_m$, be two sequences of $m$ positive integers. We denote $-\bmq=(-q_1,\ldots,-q_m)$ and $\bmp \times \bmq$ the partition which is the union of $m$ rectangles of sizes $p_i\times q_i$, namely
\[ \bmp \times \bmq = (\, \underbrace{q_1,\ldots,q_1}_{p_1 \text{ times}}\, ,\,
\underbrace{q_2,\ldots,q_2}_{p_2 \text{ times}}\, ,\,\ldots).\]
The following result was conjectured by Stanley~\cite{S1} and
proved in~\cite{F} (under a more precise statement which was also conjectured by Stanley~\cite{S2}).
\begin{theo}
For $\la=\bmp \times \bmq$ the normalized character $(-1)^k \,\widehat{\chi}^\la_{\mu,1^{n-k}}$ is a polynomial in the indeterminates $(\bmp,\bmq)$. Once $\bmq$ replaced by $-\bmq$, this polynomial has nonnegative integer coefficients.
\end{theo}

The purpose of this paper is to present a conjectured extension of this property, in the framework of Jack polynomials. 

The family of Jack polynomials $J_{\la}(\alpha)$ is indexed by partitions. It forms a basis of the algebra of symmetric functions with rational coefficients in some positive real parameter $\alpha$. We consider the transition matrix between this basis and the classical basis of power sums $p_{\rho}$. Namely we write
$$J_{\la}(\alpha)=\sum_{|\rho|= |\la|} \theta^{\la}_{\rho}(\alpha) \,p_{\rho}.$$

Let $\mu$ be a partition of weight $|\mu|=k\le |\la|=n$. Using multiplicities, we denote $\mu=(1^{m_1(\mu)},2^{m_2(\mu)},\ldots)$ and $z_\mu  = \prod_{i \ge  1} i^{m_i(\mu)} m_i(\mu)!$. We consider the quantity 
$$z_{\mu}\, \theta^{\la}_{\mu,1^{n-k}}(\alpha)$$
for which we conjecture the following positivity property.

\begin{conj}
Let $\la=\bmp \times \bmq$ and $\mu$ a partition with $m_1(\mu)=0$ and $|\mu|=k\le |\la|=n$.  
\begin{enumerate}
\item[(i)]{\setlength{\baselineskip}{1.2\baselineskip}
The quantity $z_{\mu} \theta^{\la}_{\mu,1^{n-k}}(\alpha)$ is a polynomial in the indeterminates $(\bmp,\bmq)$ and $\beta=\alpha-1$, with integer coefficients.
\par}
\item[(ii)]If $\bmq$ is replaced by $-\bmq$, the coefficients of the polynomial $(-1)^{k} \,z_{\mu}\, \theta^{\la}_{\mu,1^{n-k}}(\alpha)$ in $(\bmp,\bmq,\beta)$ are nonnegative integers. 
\item[(iii)]At least one of these coefficients is $1$.
\end{enumerate}
\end{conj}

Here two remarks are needed. Firstly, for $\alpha=1$ Jack polynomials are essentially Schur functions, and it is easily seen that
\begin{equation*}
z_{\rho}\, \theta^{\la}_{\rho}(1)=n!\,
\frac{\chi^\la_{\rho}}{\chi^\la_{1^n}},
\end{equation*}
or equivalently
\begin{equation}
\widehat{\chi}^\la_{\mu,1^{n-k}}=
\binom{n-k+m_1(\mu)}{m_1(\mu)}\,z_{\mu}\, \theta^{\la}_{\mu,1^{n-k}}(1).
\end{equation}
Therefore Stanley's ex-conjecture corresponds to the case $\alpha=1$ of our conjecture.

Secondly, expanding the Jack polynomials in terms of the ``augmented'' monomial symmetric functions, which are integral combinations of the power sums~\cite[p.110]{Ma}, and using the results of~\cite{KS2}, it is not difficult to see that $z_{\mu}\, \theta^{\la}_{\mu,1^{n-k}}(\alpha)$ is a polynomial in $\al$, hence in $\beta$. It is also easy to prove that it is a polynomial in $(\bmp,\bmq)$. 

We have checked our conjecture for $m \le 3$ and for any partition $\mu$ with $m_1(\mu)=0$ and $|\mu|-l(\mu)\le 8$. These data are available on a web page~\cite{W}. For $m=2$ after substitution of $-\bmq$ to $\bmq$, the first values are as follows :
\begin{eqnarray*}
2 \theta^{\la}_{2,1^{n-2}}(\al) & = & p_1q_1^2+p_2q_2^2+2p_1p_2q_2+p_1^2q_1+p_2^2q_2
  \\ & & \ +\beta (p_1q_1+p_2q_2+p_1q_1^2+p_2q_2^2), \\
 -3 \theta^{\la}_{3,1^{n-3}}(\al) & = & 
p_1q_1+p_2q_2+p_1q_1^3+p_2q_2^3+p_1^3q_1+p_2^3q_2\\ & & \
+3p_1^2p_2q_2+3p_1p_2^2q_2+3p_1p_2q_2^2+3p_1p_2q_1q_2+3p_1^2q_1^2+3p_2^2q_2^2\\ & & \
+\beta\, (p_1q_1+p_2q_2
+3p_1q_1^2+3p_2q_2^2+3p_1^2q_1+3p_2^2q_2+6p_1p_2q_2\\ & & \
+3p_1^2q_1^2+3p_2^2q_2^2+3p_1p_2q_2^2
+2p_1q_1^3+2p_2q_2^3+3p_1p_2q_1q_2)\\ & & \
+\beta^2\, (2p_1q_1+2p_2q_2+3p_1q_1^2+3p_2q_2^2+p_1q_1^3+p_2q_2^3).
\end{eqnarray*}

In this paper we present a general method to obtain linear identities between the coefficients $\theta_\rho^\la$. We apply this method in the simplest case $m=1$, i.e. when $\la=p \times q$, the rectangular shape formed by $p$ parts equal to $q$. 

In this situation, we prove that $(-1)^{k} \,z_{\mu}\, \theta^{\la}_{\mu,1^{n-k}}(\alpha)$ is a polynomial in $(p,-q,\beta)$, with nonnegative rational coefficients. The proof is much more cumbersome and lenghty than in the case $\alpha=1$, studied in~\cite{S1,Ra}. 

We use induction on the weight $|\mu|=k$. An explicit recurrence formula generates $(-1)^{k} \,z_{\mu}\, \theta^{\la}_{\mu,1^{n-k}}(\alpha)$ as a polynomial in $(\alpha,\beta)$, with nonnegative rational coefficients. In spite of empirical evidence, we are in lack of an argument proving that these rational numbers are actually integers. 

After substitution of $-q$ to $q$, the first values are
\begin{eqnarray*}
2 \theta^{\la}_{2,1^{pq-2}}(\al) & = & pq(\al q+p+\beta), \\
 -3 \theta^{\la}_{3,1^{pq-3}}(\al) & = & 
pq(\al q+p+\beta)(\al q+p+2\beta)+\al pq(pq+1),\\
4\theta^{\la}_{4,1^{pq-4}}(\al) & = & pq\big((\al q+p+\beta)(\al q+p+2\beta)+\al (pq+1)\big)(\al q+p+3\beta)\\ & & \
+2\al pq(pq+2)(\al q+p+\beta),\\
8\theta^{\la}_{22,1^{pq-4}}(\al) & = & 2pq(\al q+p+\beta)(\al q+p+2\beta)+2\al pq(pq+1)\\& & \
+pq(pq+2)(\al q+p+\beta)^2.
\end{eqnarray*}

We conjecture that such a property keeps true in the general case.
\begin{conj}
Let $\la=\bmp \times \bmq$ and $\mu$ a partition with $m_1(\mu)=0$ and $|\mu|=k\le |\la|=n$.
Once $\bmq$ replaced by $-\bmq$, the quantity $(-1)^{k} \,z_{\mu}\, \theta^{\la}_{\mu,1^{n-k}}(\alpha)$ has some ``natural'' expression  as a polynomial in $(\alpha,\beta)$ with nonnegative integer coefficients. 
\end{conj}

Actually our results suggest the existence of some mysterious $(\alpha,\beta)$-scheme, underlying the classical theory of Jack polynomials, where $\beta$ would play a role as important as $\al$. 
However we have no conjectured expression, nor any combinatorial interpretation, giving the quantity $ (-1)^{k} \,z_{\mu}\, \theta^{\la}_{\mu,1^{n-k}}(\alpha)$  as a polynomial in $(\alpha,\beta)$, with nonnegative integer coefficients. 

The paper is organized as follows. In Section 2 we introduce our notations and recall general facts about (shifted) symmetric functions and (shifted) Jack polynomials. In Section 3, starting from the generalized binomial formula, we define an isomorphism between symmetric and shifted symmetric functions. In Sections 4 and 5 we use this method to obtain several linear identities between the $\theta_\rho^\la$. Section 6 is devoted to the case $m=1$. In Section 7 we compare our proof with those previously given for $\al=1$~\cite{S1,Ra}. Finally a conjectural $(\alpha,\beta)$-extension of Jack polynomials is discussed in Section 8.

\section{(Shifted) Jack polynomials}

The standard reference for symmetric functions and Jack polynomials are~\cite[Section 6.10]{Ma} and~\cite{S0}. Although the theory of symmetric functions goes back to the early 19th century, the 
notion of ``shifted symmetric'' functions is quite recent. We refer 
to~\cite{KS,Ok1,Ok2,Ok3} and to other references given there.

\subsection{Symmetric functions}

Let $x=\{x_1,x_2,x_3,\ldots\}$ be an infinite set of indeterminates, and $\mathcal{S}$ the corresponding algebra of symmetric functions
with coefficients in $\mathbf{Q}$.
Let $\mathbf{Q}[\alpha]$ be the field of rational functions in
some indeterminate $\alpha$ (wich may be considered as a positive real number), and $\mathbf{S}=\mathcal{S}\otimes\mathbf{Q}[\alpha]$
the algebra of symmetric functions with coefficients in $\mathbf{Q}[\alpha]$. The parameter $\al$ being kept fixed, for clarity of display, we shall omit its dependence in any notation below.

A partition $\la= (\la_1,...,\la_n)$
is a finite weakly decreasing
sequence of nonnegative integers, called parts. The number
$l(\la)$ of positive parts is called the length of
$\la$, and $|\la| = \sum_{i = 1}^{n} \la_i$
the weight of $\la$. For any integer $i\geq1$,
$m_i(\la) = \textrm{card} \{j: \la_j  = i\}$
is the multiplicity of the part $i$ in $\la$.  Clearly
$l(\la)=\sum_{i\ge1} m_i(\la)$ and
$|\la|=\sum_{i\ge1} im_i(\la)$. We also write
$\la= (1^{m_1(\la)},2^{m_2(\la)},3^{m_3(\la)},\ldots)$ and set
\[z_\la  = \prod_{i \ge  1} i^{m_i(\lambda)} m_i(\lambda) ! .\]
Being given two partitions, we write $\mu \subseteq \la$ if $\mu_i \le \la_i$ for any $i$. We denote $\la^{'}$ the partition conjugate to $\la$, whose parts are given by $m_i(\la^{'})=\la_i-\la_{i+1}$. We
have $\la^{'}_i=\sum_{j \ge i} m_j(\la)$. 

We define
\begin{equation*}
\begin{split}
h_\la&=\prod_{(i,j) \in \la} \left(\la^{'}_j-i+1+\al(\la_i-j)\right),\\
h^{'}_\la&=\prod_{(i,j) \in \la}  \left(\la^{'}_j-i+\al(\la_i-j+1)\right),\\
(u)_{\la}&=\prod_{(i,j) \in \la} \left(u+j-1-(i-1)/\al\right).
\end{split}
\end{equation*}
The last quantity is a generalization of the ``raising'' factorial, in terms of the ``$\al$-contents'' $j-1-(i-1)/\al$.

The power sum symmetric functions
are defined by $p_{k}(x)=\sum_{i \ge 1} x_i^k$.
They form an algebraic basis of $\mathbf{S}$. A linear basis is given by the symmetric functions
\[p_{\la}=\prod_{i=1}^{l(\la)}p_{\la_{i}}=\prod_{i\geq1}p_i{}^{m_{i}(\la)}.\]
The algebra $\mathbf{S}$ may be endowed with a scalar product 
$<\, , \,>$ for which we have two orthogonal bases, both indexed by partitions :
\begin{enumerate}
\item[(i)] the basis of power sum symmetric functions, with
\[<p_\la,p_\mu>=\delta_{\la \mu}\,\alpha^{l(\la)} z_{\la},\]
\item[(ii)] the basis of Jack symmetric functions, with
\[<J_\la,J_\mu>=\delta_{\la \mu}\,h_\la h^{'}_\la.\]
\end{enumerate}
We write $\theta^{\la}_{\rho}$ for the transition matrix between these two orthogonal bases, namely
$$J_{\la}=\sum_{|\rho|= |\la|} \theta^{\la}_{\rho} \,p_{\rho}.$$

If we restrict to a finite set of $N$
indeterminates $x=(x_1,\ldots,x_N)$, we have
\begin{equation}
J_\la(1^N):=J_\la(1,\ldots,1)=\al^{|\la|} \,(N/\al)_\la.
\end{equation}
Denoting $j_{\la}=h_\la h^{'}_\la$, we introduce
\[J_{\la}^{\sharp}=\frac{J_\la}{j_\la}, \quad \quad
J_\la^{\star}=\frac{J_\la}{J_\la(1^N)}.\]
The first relation defines the basis dual of $J_{\la}$ with respect to the scalar product $<\, , \,>$. In contrast with $J_\la^{\star}$, it does not depend of $N$. We have
\begin{equation*}
\prod_{i,j=1}^N (1-x_iy_j)^{-1/\al}
=\sum_{\la} J_{\la}(x)\, J_{\la}^{\sharp}(y).
\end{equation*}

\subsection{Shifted symmetric functions}

A polynomial in $N$ indeterminates $x=(x_1,\ldots,x_N)$ with coefficients in $\mathbf{Q}[\alpha]$ is said to be ``shifted symmetric'' if it is symmetric in the $N$ ``shifted variables'' 
$x_i-i/\alpha$. 

Dealing with an infinite set of indeterminates $x=\{x_1,x_2,x_3,\ldots\}$, in analogy with symmetric functions, 
a ``shifted symmetric function'' $f$ is a family $\{f_i, i\ge 1\}$ with the two following properties :
\begin{enumerate}
    \item[(i)]  $f_i$ is shifted symmetric in $(x_1,x_2,\ldots,x_i)$,
    \item[(ii)]  $f_{i+1}(x_1,x_2,\ldots,x_i,0)=f_i(x_1,x_2,\ldots,x_i)$.
\end{enumerate}
This defines $\mathbf{S}^{\ast}$, the algebra of shifted symmetric functions with coefficients in $\mathbf{Q}[\alpha]$. A typical example is provided by the ``shifted power sums''
\[p_k^{\star}(x)=\sum_{i\ge 
1}\Big([x_i-(i-1)/\alpha]_k-[-(i-1)/\alpha]_k\Big),\]
with $[x]_k=x(x-1)\cdots(x-k+1)$. These shifted symmetric functions generate ${\mathbf{S}}^{\ast}$ algebraically.

Any element $f\in \mathbf{S}^{\ast}$ may be evaluated at any sequence 
$x=(x_1,x_2,\ldots)$ with finitely many non zero terms, hence at any partition $\la$. Moreover by analyticity, $f$ is entirely determined by its restriction $f(\la)$ to partitions. This identification is usually performed and $\mathbf{S}^{\ast}$ is considered as a function 
algebra on the set of partitions.

For any partition $\mu$ there exists a shifted symmetric function 
$J_\mu^{\dag}$ such that
\begin{enumerate}
\item[(i)]  degree $J_\mu^{\dag} = |\mu|$,
\item[(ii)] $J_\mu^{\dag}(\la)=0$ except if 
$\mu \subseteq \la$, and 
$J_\mu^{\dag}(\mu) \neq 0$.
\end{enumerate}
It is a very remarkable fact that in this definition, 
the overdetermined system of linear conditions (ii) may be replaced by
the weaker conditions
\begin{enumerate}
\item[(iii)] $J_\mu^{\dag}(\la)=0$ except if $|\mu| \le |\la|$, and $J_\mu^{\dag}(\mu) \neq 0$.
\end{enumerate}
The function $J_\mu^{\dag}$ is called the ``shifted Jack polynomial'' associated with $\mu$. It is unique up to the value of $J_\mu^\dag(\mu)$. 

A map $\mathbf{S}^{\ast} \rightarrow \mathbf{S}$ can be defined, which associates to any shifted symmetric function $f\in \mathbf{S}^{\ast}$ its ``leading symmetric term'' denoted $[f]$. By definition $[f]$ is the highest degree term of $f$, which is necessarily symmetric.

It is another very remarkable fact that Jack polynomials are the leading symmetric terms of shifted Jack polynomials. More precisely we have
\[\left[\frac{J_\mu^\dag}{J_\mu^\dag(\mu)}\right]= \al^{|\mu|}J^{\sharp}_\mu.\]
Hence the family $\{J_\mu^{\dag}/J_\mu^\dag(\mu)\}$ forms a basis of the algebra $\mathbf{S}^\ast$.

\subsection{Generalized binomial formula}

Jack polynomials allow to write the following generalization of the classical binomial formula
\begin{equation}
J^{\star}_\la(1+x_1,\ldots,1+x_N)=\sum_{\mu \subseteq \la} \binom{\la}{\mu} J^{\star}_\mu (x_1,\ldots,x_N),
\end{equation}
which was first studied independently in~\cite{Ka,La1}.
\newpage
The generalized binomial coefficients thus introduced may be given the following alternative definition
\begin{equation}
\mathrm{exp}(p_1) J_{\mu}^{\sharp} =\sum_{\la \supseteq \mu} \al^{|\la|-|\mu|} \binom{\la}{\mu} J_{\la}^{\sharp}.
\end{equation}
The equivalence of both properties was proved in~\cite{La2}, as the limit of a more general result, obtained in the framework of Macdonald polynomials. This second definition has the advantage of being independent of $N$.

It was first observed in~\cite{Ok1} that the 
generalized binomial coefficient $\binom{\la}{\mu}$ is merely the shifted Jack polynomial $J_\mu^\dag(\la)$ suitably normalized:
\begin{equation}
\frac{J_\mu^\dag(\la)}{J_\mu^\dag(\mu)} = \binom{\la}{\mu}.
\end{equation}
This property is actually a special case of a more general correspondence, that we shall explicitate in the next section.

\section{Symmetric versus shifted symmetric}

Being given any symmetric function $f \in \mathbf{S}$, we write
\begin{equation}
\mathrm{exp}(p_1) f =\sum_{\la} \al^{|\la|} f^\#(\la) J_{\la}^{\sharp},
\end{equation}
i.e. we develop the inhomogeneous symmetric left-hand side in terms of the Jack polynomials basis.

Then the results recalled in Section 2 can be rephrased as follows.
\begin{theo}
The coefficient $f^\#(\la)$ is a shifted symmetric function of $\la$. The map $f \rightarrow f^\#$ is an isomorphism of $\mathbf{S}$ onto $\mathbf{S}^{\ast}$. If $f$ is homogeneous, one has $[f^\#]=f$.
\end{theo}
A direct proof would be possible, but it is out of the scope of this paper. Here we shall only mention that relations (2.3) and (2.4) imply
\[\left(\al^{|\mu|}J^{\sharp}_\mu\right)^\#=
\frac{J_\mu^\dag}{J_\mu^\dag(\mu)},\]
from which follows
\[f^\#=\sum_\mu \al^{-|\mu|} <f,J^{\sharp}_\mu>
\frac{J_\mu^\dag}{J_\mu^\dag(\mu)}.\]
Hence the statements. Observe that if $f$ is not homogeneous, $[f^\#]$ is its highest degree term.

We now give some examples and properties of this isomorphism.
\begin{prop}
Let $f \in \mathbf{S}$ be a symmetric function, homogeneous of degree $k$. For any positive integer $r$, we have
\begin{equation*}
\begin{split}
\mathrm{(i)}& \quad  \left(\frac{p_1^r}{r!}f\right)^\#(\la)= 
\binom{|\la|-k}{r}f^\#(\la),\\
\mathrm{(ii)}& \quad \binom{|\la|-k}{r-k}f^\#(\la)=
\sum_{|\rho|=r}\binom{\la}{\rho} f^\#(\rho).
\end{split}
\end{equation*} 
\end{prop}
\begin{proof}
By the definition (3.1) for any $s \ge 0$ we have
\[\frac{p_1^s}{s!} f=\sum_{|\la|=k+s} \al^{|\la|} f^\#(\la) J_{\la}^{\sharp}.\]
Hence
\[\mathrm{exp}(p_1) \frac{p_1^r}{r!} f =
\sum_{s\ge 0} \sum_{|\la|=k+r+s} \frac{(r+s)!}{r!s!} 
 \al^{|\la|} f^\#(\la) J_{\la}^{\sharp}.\]
Hence (i). On the other hand, we have
\begin{equation*}
\begin{split}
\mathrm{exp}(p_1) \frac{p_1^r}{r!} f=&
\sum_{\rho} \al^{|\rho|} f^\#(\rho)\, \frac{p_1^r}{r!}\, J_{\rho}^{\sharp}\\
=&\sum_{\rho} \al^{|\rho|} f^\#(\rho) \sum_{|\la|=|\rho|+r} \al^r
 \binom{\la}{\rho} J_{\la}^{\sharp},
\end{split}
\end{equation*}
where the second equality is a straightforward consequence of (2.3). In other words,
\[\left(\frac{p_1^r}{r!}f\right)^\#(\la)=
\sum_{|\rho|=|\la|-r} \binom{\la}{\rho} f^\#(\rho).\]
Comparing with (i), we obtain (ii).
\end{proof}

In particular writing (ii) with $f=\al^{|\mu|}J^{\sharp}_\mu$ yields
\[\binom{|\la|-|\mu|}{r-|\mu|}\binom{\la}{\mu}=
\sum_{|\rho|=r}\binom{\la}{\rho} \binom{\rho}{\mu}.\]
We are interested in the isomorphism $f \rightarrow f^\#$ because of the following important example.
\begin{prop}
Let $\mu$ be a partition with weight $|\mu|=k$. For any partition $\la$ with $|\la|=n \ge k$, we have
\begin{equation}
\left(\al^{k-l(\mu)} p_\mu\right)^\#(\la)=
\binom{n-k+m_1(\mu)}{m_1(\mu)} \,z_{\mu}\, \theta^{\la}_{\mu,1^{n-k}}.
\end{equation}
Thus $\theta^{\la}_{\mu,1^{n-k}}$ is a shifted symmetric function of $\la$.
\end{prop}
\begin{proof}
For any partition $\rho$, by orthogonality of the power sums we have
\[<J_\la,p_{\rho}>=
\al^{l(\rho)}\,z_{\rho}\,\theta^{\la}_{\rho}.\]
Hence by orthogonality of the Jack polynomials,
\[p_{\rho}=\al^{l(\rho)}\,z_{\rho} \sum_{|\la|=n} \theta^{\la}_{\rho}\, J_{\la}^{\sharp}. \]
By the definition (3.1) for $\rho=(\mu,1^{n-k})$ this implies
\[(p_{\mu}){}^{\#}(\la)=
\al^{l(\mu)-k} \frac{z_{\mu,1^{n-k}}}{(n-k)!}\,\theta^{\la}_{\mu,1^{n-k}}.\]
\end{proof}
\begin{coro}
Let $\mu$ be a partition with weight $|\mu|=k$. For any partition $\la$ with $|\la|=n \ge k$ and any $r\ge 0$, we have
\begin{equation}
\binom{n-k+m_1(\mu)}{n-r} \theta^{\la}_{\mu,1^{n-k}}=
\sum_{|\rho|=r}\binom{\la}{\rho} \theta^{\rho}_{\mu,1^{r-k}}.
\end{equation}
\end{coro}
\begin{proof}
A consequence of Proposition 1 (ii).
\end{proof}

Taking Proposition 2 into account, we may state an alternative formulation of our conjecture, which gets rid of the restriction $m_1(\mu)=0$. The quantity $(\al^{k-l(\mu)} p_\mu){}^\#(\bmp \times \bmq)$ being a shifted symmetric function of $\bmp \times \bmq$, it is a polynomial in the indeterminates $(\bmp,\bmq)$ with coefficients in $\mathbf{Q}[\alpha]$ (as mentioned in the introduction, they are actually polynomials in $\al$).
\begin{konj}
Let $\la=\bmp \times \bmq$ and $\mu$ a partition with $|\mu|=k\le |\la|=n$.  
\begin{enumerate}
\item[(i)]The shifted symmetric function  $(\al^{k-l(\mu)} p_\mu){}^\#(\bmp \times \bmq)$ is a polynomial in $(\bmp,\bmq,\beta)$ with integer coefficients.
\item[(ii)]If $\bmq$ is replaced by $-\bmq$, the coefficients of the polynomial $(-1)^{k} \,(\al^{k-l(\mu)} p_\mu){}^\#(\bmp \times \bmq)$ in $(\bmp,\bmq,\beta)$ are nonnegative integers. 
\item[(iii)]At least one of these coefficients is $1$.
\end{enumerate}
\end{konj}
This formulation is closer to the original Stanley's conjecture. Actually for $\al=1$, equation (1.1) can be written
\[\widehat{\chi}^\la_{\mu,1^{n-k}}= (p_{\mu}){}^\#(\la).\]
Indeed for $\al=1$ we have $J_\la= h_\la s_\la$, with $s_\la$ the Schur function, and the result follows from the classical Frobenius formula and the dimension formula $\chi^\la_{1^n}=n!/h_\la$~\cite[p.114 and 117]{Ma}.

\section{Linear identities}

\subsection{A method}

Let $D$ be some linear operator (not necessarily differential) acting on $\mathbf{S}$. Let us assume that the action of $D$ on Jack polynomials is explicitly known, i.e.
\[D J_\rho^{\sharp}=
\sum_{\sigma} a_{\rho\sigma}J_\sigma^{\sharp}.\]
On the one hand, we may apply $D$ to (3.1), which yields 
\begin{equation*}
D(\mathrm{exp}(p_1) f) =\sum_{\la} \al^{|\la|} f^\#(\la) D J_{\la}^{\sharp}=\sum_{\la,\rho}  \al^{|\la|} f^\#(\la) \,a_{\la\rho}J_\rho^{\sharp}.
\end{equation*}
On the other hand $D(\mathrm{exp}(p_1) f)$ may be expressed by using the Baker-Campbell-Hausdorff formula
\[D(\mathrm{exp}(p_1) f)=\mathrm{exp}(p_1) \left(Df+[D,p_1]f+\frac{1}{2!}[[D,p_1],p_1]f +\frac{1}{3!}[[[D,p_1],p_1],p_1]f 
+\ldots\right).\]

By comparison we obtain the following linear identity between shifted symmetric functions:
\begin{equation*}
\Big(Df+[D,p_1]f+\frac{1}{2!}[[D,p_1],p_1]f +\frac{1}{3!}[[[D,p_1],p_1],p_1]f 
+\ldots\Big)^\#(\la)=\sum_{\rho}\al^{|\rho|-|\la|}a_{\rho\la}f^\#(\rho).
\end{equation*}

If moreover the action of $D$ is explicitly known on the power sums, by specializing $f= \al^{|\mu|-l(\mu)} p_\mu$ and using (3.2), we shall obtain a linear identity between some coefficients $\theta^{\la}_{\rho}$. Many examples are given below.

\subsection{Pieri formula}

For any partition $\la$ and any integer $1 \le i \le l(\la)+1$, we denote $\la^{(i)}$ the partition $\mu$ (if it exists) such that $\mu_j=\la_j$ for $j\neq i$ and $\mu_i=\la_i +1$. Similarly for any integer $1 \le i \le l(\la)$, we denote $\la_{(i)}$ the partition $\nu$ (if it exists) such that $\nu_j=\la_j$ for $j\neq i$ and $\nu_i=\la_i -1$. 

Jack polynomials satisfy the following generalization of Pieri formula~\cite{Ma,S0}: 
\[p_1 \,J_\la =\sum_{i=1}^{l(\la)+1} 
c_i(\la) \, J_{\la^{(i)}}.\]
The Pieri coefficients $c_i(\la)$ have the following analytic 
expression~\cite{La0}, see also~\cite[Prop.5, p.299]{La2}:
\[c_i(\la) = \frac {1}{\alpha \la_i+l(\la)-i+2}
\prod_{\begin{subarray}{c}j=1 \\ j \neq i\end{subarray}}^{l(\la)+1} 
\frac{\alpha(\la_i-\la_j)+j-i+1}
{\alpha(\la_i-\la_j)+j-i}.\]
In~\cite[p.300]{La2} (see also \cite{La1}) we have proved that these coefficients are connected with the generalized binomial coefficients $\binom{\la^{(i)}}{\la}$ by
\begin{equation}
c_i(\la)=\al \binom{\la^{(i)}}{\la}\, \frac{j_{\la}}{j_{\la^{(i)}}} . \end{equation}
From which follows
\[\binom{\la}{\la_{(i)}} = \left(\la_i +\frac {l(\la)-i}{\alpha}\right)
\prod_{\begin{subarray}{c}j=1 \\ j \neq i \end{subarray}}^{l(\la)} 
\frac{\alpha(\la_i-\la_j)+j-i-1}
{\alpha(\la_i-\la_j)+j-i}.\]

The method of Section 4.1 is very easy to apply for $D=p_1$ or $D=\partial/\partial p_1$, since in these cases we have $[D,p_1]=0$ or $1$, respectively. We write $D^\bot$ for the adjoint of any linear operator $D$ with respect to the scalar product $<\, , \,>$.
\begin{prop}
For any symmetric function $f \in \mathbf{S}$ we have
\begin{equation*}
\begin{split}
\mathrm{(i)}& \quad
(p_1f)^\#(\la) =
\sum_{i=1}^{l(\la)} \binom{\la}{\la_{(i)}} \, f^\#(\la_{(i)}),\\
\mathrm{(ii)}& \quad
\left(\frac{\partial}{\partial p_1} f\right)^\#(\la) +f^\#(\la)=
\sum_{i=1}^{l(\la)+1} c_i(\la) \, f^\#(\la^{(i)}).
\end{split}
\end{equation*}
If $f$ is homogeneous of degree $k$, $\mathrm{(i)}$ is equal to $(|\la|-k)f^\#(\la)$.
\end{prop}
\begin{proof}
To prove (i) we may use Proposition 1 or apply $p_1$ to (3.1), which yields
\begin{equation*}
\begin{split}
p_1\mathrm{exp}(p_1) f&=\sum_{\la} \al^{|\la|} f^\#(\la)
\sum_{i=1}^{l(\la)+1} c_i(\la) \,
\frac{j_{\la^{(i)}}}{j_{\la}}\, J_{\la^{(i)}}^{\sharp}\\
&=\sum_{\la} \al^{|\la|} J_{\la}^{\sharp}
\sum_{i=1}^{l(\la)} \binom{\la}{\la_{(i)}} \, f^\#(\la_{(i)}),
\end{split}
\end{equation*}
where the second equality is a consequence of (4.1).

To prove (ii)  observe that $p_1^\bot=\al \partial/\partial p_1$. This can be proved by the same argument than in~\cite[p.76]{Ma}: by linearity, it is enough to check the fact on power sums. But in~\cite[Prop.11, p.306]{La2} we have shown that
\[p_1^\bot J_\la= \al \sum_{i=1}^{l(\la)} \binom{\la}{\la_{(i)}} \,J_{\la_{(i)}}.\]
Applying $p_1^\bot$ to (3.1) we thus obtain
\begin{equation*}
\begin{split}
\frac{\partial}{\partial p_1}(\mathrm{exp}(p_1) f)
&=\mathrm{exp}(p_1) \frac{\partial}{\partial p_1}f+\mathrm{exp}(p_1) f\\
&=\sum_{\la} \al^{|\la|} f^\#(\la)
\sum_{i=1}^{l(\la)} \binom{\la}{\la_{(i)}} \,
\frac{j_{\la_{(i)}}}{j_{\la}}\, J_{\la_{(i)}}^{\sharp}\\
&=\sum_{\la} \al^{|\la|} J_{\la}^{\sharp}
\sum_{i=1}^{l(\la)+1} c_i(\la) \, f^\#(\la^{(i)}).
\end{split}
\end{equation*}
\end{proof}

A similar argument may be used for the differential operator
$E_0=\sum_{i=1}^N \partial/\partial x_i$,
which is dependent of the number of variables $N$. 
\begin{prop}
For any symmetric function $f \in \mathbf{S}$ we have
\begin{equation*}
(E_0f)^\#(\la) +N f^\#(\la)=
\sum_{i=1}^{l(\la)+1} c_i(\la) \,(N+\al\la_i-i+1)\, f^\#(\la^{(i)}).
\end{equation*}
\end{prop}
\begin{proof}Applying  $E_0$ to the generalized binomial formula (2.3), we obtain easily
\[E_0 J_\la^\star= \sum_{i=1}^{l(\la)} \binom{\la}{\la_{(i)}} \,J^\star_{\la_{(i)}}.\]
Hence
\begin{equation*}
\begin{split}
E_0(\mathrm{exp}(p_1) f)&=\mathrm{exp}(p_1) E_0f+N\mathrm{exp}(p_1) f
\\
&=\sum_{\la} \al^{|\la|} f^\#(\la) \sum_{i=1}^{l(\la)} \binom{\la}{\la_{(i)}} \,\frac{j_{\la_{(i)}}}{j_{\la}}\,
\frac{J_\la(1^N)}{J_{\la_{(i)}}(1^N)}\,J^\sharp_{\la_{(i)}}\\
&=\sum_{\la} \al^{|\la|} J_{\la}^{\sharp}
\sum_{i=1}^{l(\la)+1} c_i(\la) \, \frac{J_{\la^{(i)}}(1^N)}{J_\la(1^N)} f^\#(\la^{(i)}).
\end{split}
\end{equation*}
By (2.1) we have $J_{\la^{(i)}}(1^N)/J_\la(1^N)=N+\al\la_i-i+1$. 
\end{proof}

\subsection{Other examples}

For any $k\ge 0$ we introduce the differential operators
\begin{equation*}
\begin{split}
E_k&=\sum_{i=1}^N x_i^k\, \frac{\partial}{\partial x_i},\\
D_k&=\sum_{i=1}^N x_i^k\, \frac{\partial^2}{\partial x_i^2}
+\frac{2}{\al} \sum_{\begin{subarray}{c}i,j=1\\i\neq j\end{subarray}}^N 
\frac{x_i^k}{x_i-x_j}\frac{\partial}{\partial x_i}.
\end{split}
\end{equation*} 
For $k\neq 0$, $E_k$ is independent of $N$. It is not difficult to check
\begin{equation*}
\begin{split}
D_1&=\frac{1}{2}[E_0,D_2],\qquad \quad D_0=[E_0,D_1],\\
E_2&= \sum_{k \ge 1} k\, p_{k+1}\, \frac{\partial}{\partial p_k}=
\frac{1}{2} [D_2,p_1]-\frac{1}{\al}(N-1)p_1,\\
E_2^\bot&= \sum_{k \ge 1} (k+1) p_{k}\, \frac{\partial}{\partial p_{k+1}}=E_0-N\frac{\partial}{\partial p_1}.
\end{split}
\end{equation*}
If we write 
\[d_1(\la)=\sum_{(i,j) \in \la} \big(j-1-(i-1)/\al\big)\]
for the sum of the $\alpha$-contents of $\la$, it is well known~\cite[p.84]{S0} that the Jack polynomials are eigenfunctions of $D_2$, namely
\[D_2 J_\la=2\big(d_1(\la)+|\la|(N-1)/\al\big)J_\la.\]
This implies
\begin{equation*}
\begin{split}
D_1 J_\la^\star&= \sum_{i=1}^{l(\la)} \binom{\la}{\la_{(i)}} \,\big(\la_i-1+(N-i)/\al\big)\, J^\star_{\la_{(i)}},\\
E_2 J_\la&=\sum_{i=1}^{l(\la)+1} 
c_i(\la) \,\big(\la_i-(i-1)/\al\big)\, J_{\la^{(i)}},\\
E_2^\bot J_\la&= \al \sum_{i=1}^{l(\la)} \binom{\la}{\la_{(i)}} \,\big(\la_i-1-(i-1)/\al\big)\, J_{\la_{(i)}}.
\end{split}
\end{equation*} 

We may apply the method of Section 4.1 to $E_2$ and $D_1$. 
\begin{prop}
For any symmetric function $f \in \mathbf{S}$ we have
\begin{equation*}
(E_2f)^\#(\la) +(p_2f)^\#(\la)=
\sum_{i=1}^{l(\la)} \binom{\la}{\la_{(i)}} \,
\big(\la_i-1-(i-1)/\al\big)\, f^\#(\la_{(i)}).
\end{equation*}
\end{prop}
\begin{proof}
Using $[E_2,p_1]=p_2$, we have
\begin{equation*}
\begin{split}
E_2(\mathrm{exp}(p_1) f)&= \mathrm{exp}(p_1) (E_2f+p_2f)\\
&=\sum_{\la} \al^{|\la|} f^\#(\la)
\sum_{i=1}^{l(\la)+1} c_i(\la)  \,\big(\la_i-(i-1)/\al\big)\, \frac{j_{\la^{(i)}}}{j_\la}J_{\la^{(i)}}^\sharp\\
&=\sum_{\la} \al^{|\la|} J_{\la}^{\sharp}
\sum_{i=1}^{l(\la)} \binom{\la}{\la_{(i)}} \,\big(\la_i-1-(i-1)/\al\big)\, f^\#(\la_{(i)}).
\end{split}
\end{equation*}
\end{proof}

\begin{prop}
For any symmetric function $f \in \mathbf{S}$ we have
\begin{multline*}
\Big(\big(D_1+2E_1+p_1 +N(N-1)/\al\big)f\Big)^\#(\la)=\\ \al
\sum_{i=1}^{l(\la)+1} c_i(\la) \,\big(\la_i+(N-i)/\al\big)\,
\big(\la_i+(N-i+1)/\al\big)\, f^\#(\la^{(i)}).
\end{multline*}
\end{prop}
\begin{proof}
We use $[D_1,p_1]=2 E_1+N(N-1)/\al$,
hence $[[D_1,p_1],p_1]=2p_1$. We have
\begin{equation*}
\begin{split}
D_1(\mathrm{exp}(p_1) f)&= \mathrm{exp}(p_1) (D_1f+2E_1f+p_1f+\frac{1}{\al}N(N-1)f)\\
&=\sum_{\la} \al^{|\la|} f^\#(\la) \sum_{i=1}^{l(\la)}
\binom{\la}{\la_{(i)}} \,\big(\la_i-1+(N-i)/\al\big)\, 
\frac{J_\la(1^N)}{J_{\la_{(i)}}(1^N)}\,
\frac{j_{\la_{(i)}}}{j_\la}\,
J^\sharp_{\la_{(i)}}\\
&=\sum_{\la} \al^{|\la|} J_{\la}^{\sharp}
\sum_{i=1}^{l(\la)+1} c_i(\la) \,\big(\la_i+(N-i)/\al\big)\,
\frac{J_{\la^{(i)}}(1^N)}{J_\la(1^N)}\, f^\#(\la^{(i)}).
\end{split}
\end{equation*}
\end{proof}

\subsection{Generalization}

The previous results may be extended as follows. Let $\Delta_0= p_1$ , $\nabla_0=E_0$ and for any $k > 0$ define
\begin{equation*}
\Delta_k= \frac{1}{2}[D_2,\Delta_{k-1}],\quad \quad
\nabla_k= \frac{1}{2}[\nabla_{k-1},D_2].
\end{equation*} 
Clearly we have 
\[\Delta_1=E_2+(N-1)p_1/\al, \quad \quad \nabla_1=D_1.\]
Both families are in correspondence through $D\rightarrow D^\bot$. Starting from
\[E_0^\bot=E_2+Np_1/\al,\] 
we readily obtain
\[\nabla_k^\bot=\Delta_{k+1}+\Delta_{k}/\al.\] 
The operators $\Delta_k$ and $\nabla_k$ are respectively increasing and decreasing degree by 1. We have
\begin{equation*}
\begin{split}
\Delta_k J_\la&=\sum_{i=1}^{l(\la)+1} 
c_i(\la) \,\big(\la_i+(N-i)/\al\big)^k\, J_{\la^{(i)}},\\
\nabla_k J_\la^\star&= \sum_{i=1}^{l(\la)} \binom{\la}{\la_{(i)}} \,\big(\la_i-1+(N-i)/\al\big)^k\, J^\star_{\la_{(i)}},
\end{split}
\end{equation*} 
which are easy consequences of $d_1(\la^{(i)})-d_1(\la)=\la_i-(i-1)/\al$.

In this section we have applied our method to the operators $\Delta_k$ and $\nabla_k$ with $k=0,1$. The case $k=2$ is also easy to handle. For instance we have
\begin{equation*}
\begin{split}
&\Delta_2=D_3+E_2+(N-1)^2p_1/\al^2,\\
[&\Delta_2,p_1]=2E_3+(2N-3)p_2/\al+p_1^2/\al +p_2,\\
[[&\Delta_2,p_1],p_1]=2p_3.
\end{split}
\end{equation*} 
For higher values of $k$, the computations become very messy.

The operator $D_2$ is degree preserving and must be treated separately. Defining
\[D_2^\dag = D_2-2(N-1) E_1/\al,\]
we have obviously $D_2^\dag J_\la=2 d_1(\la) J_\la$. Hence $D_2^\dag$ is independent of $N$.
\begin{prop}
For any symmetric function $f \in \mathbf{S}$ we have
\[\left(\big(D_2^\dag+2E_2+p_2 \big)f\right)^\#(\la)
=2 d_1(\la)\,f^\#(\la).\]
\end{prop}
\begin{proof}
Consequence of $[D_2^\dag,p_1]=2 E_2$ and $[D_2^\dag,p_1],p_1]=2p_2$.
\end{proof}

For completeness, we mention that using $p_2=[E_2,p_1]$ and $D_0=[E_0,D_1]$, the operators $p_2$, $D_0$, and $p_2^\bot=\al^2 \partial/\partial p_2$ (which are respectively increasing and decreasing degree by 2) can also be managed very easily. This is left to the reader.

\section{Specialization to power sums}

We now specialize the results of Section 4 for $f=p_\mu$ with $|\mu|=k$. From now on, without any loss of generality, and to avoid superfluous complications, we shall assume $m_1(\mu)=0$. In order to simplify notations, we shall also write 
\[\vartheta^{\la}_{\mu}= z_{\mu}\, \theta^{\la}_{\mu,1^{n-k}}.\]
Then Proposition 2 reads $(p_\mu)^\#(\la)=\al^{l(\mu)-k}\, \vartheta^{\la}_{\mu}$.

\subsection{A lemma}

The following result is proved by an easy induction on $N$.

\begin{lem} For any integer $r\ge 2$, we have
\[2\sum_{\begin{subarray}{c}i,j=1\\i\neq j\end{subarray}}^N 
\frac{x_i^r}{x_i-x_j}=\sum_{i=1}^{r-2}p_i\,p_{r-i-1}+(2N-r)p_{r-1}.\]
\end{lem}

By straightforward computation we obtain a corollary which will be central for our purpose.
\begin{prop} 
Denote $\beta=\al-1$. For any integer $j\ge 0$, we have
\begin{multline*}
D_jp_\mu=p_\mu \Big(
\sum_{r, s} rsm_r(\mu)(m_s(\mu)-\delta_{rs}) 
\frac{p_{r+s+j-2}}{p_rp_s}
+\frac{2N-j}{\al}\sum_{r}rm_r(\mu)\frac{p_{r+j-2}}{p_r}\\+
\frac{\beta}{\al} \sum_{r} r(r-1)m_r(\mu) \frac{p_{r+j-2}}{p_r}
+\frac{1}{\al}\sum_r rm_r(\mu) \sum_{i=1}^{r+j-3} 
\frac{p_ip_{r-i+j-2}}{p_r}\Big).\quad
\end{multline*}
\end{prop}

For clarity of display, the following notations will be useful. For any integer $r\ge 2$ we denote $\mu_{\downarrow (r)}$ the partition (if it exists) obtained from $\mu$ by substracting one part $r$ and adding one part $r-1$. Similarly we denote $\mu_{\uparrow (r)}$ the partition (if it exists) obtained from $\mu$ by substracting one part $r$ and adding one part $r+1$. We have $|\mu_{\downarrow (r)}|=k-1$, $|\mu_{\uparrow (r)}|=k+1$, $l(\mu_{\downarrow (r)})=l(\mu_{\uparrow (r)})=l(\mu)$.

For any integers $r,s\ge 1$ we denote 
$\mu_{\downarrow (rs)}$ the partition (if it exists) obtained from $\mu$ by substracting one part $r$ and one part $s$, and adding one part $r+s-1$. Similarly we write $\mu_{\uparrow (rs)}$ the partition (if it exists) obtained from $\mu$ by adding one part $r$ and one part $s$, and substracting one part $r+s+1$. We have $|\mu_{\downarrow (rs)}|=|\mu_{\uparrow (rs)}|=k-1$, $l(\mu_{\downarrow (rs)})=l(\mu)-1$ and $l(\mu_{\uparrow (rs)})=l(\mu)+1$.

Writing Proposition 8 for $j=1$ we readily obtain
\begin{multline}
D_1p_\mu=\sum_{r, s} rsm_r(\mu)(m_s(\mu)-\delta_{rs}) p_{\mu_{\downarrow (rs)}}
+\frac{2N-1}{\al}\sum_{r}rm_r(\mu)p_{\mu_{\downarrow (r)}}\\  
\qquad \qquad +
\frac{\beta}{\al}\sum_{r} r(r-1)m_r(\mu) p_{\mu_{\downarrow (r)}}+
\frac{1}{\al}\sum_r rm_r(\mu)\sum_{i=1}^{r-2}p_{\mu_{\uparrow (i,r-i-1)}}.
\end{multline}

\subsection{Raising degree}
We first consider the cases of operators raising degree, i.e. $p_1$ and $E_2$. Proposition 3 (i) writes
\begin{equation}
(n-k) \vartheta^{\la}_{\mu} =
\sum_{i=1}^{l(\la)} \binom{\la}{\la_{(i)}} \, \vartheta^{\la_{(i)}}_{\mu}.
\end{equation}
Observe that it is only (3.3) written for $r=n-1$.

An easy computation gives
\[E_2p_\mu=\sum_{r} rm_r(\mu)\,p_{\mu_{\uparrow (r)}},\]
so that Proposition 5 yields
\begin{equation}
\vartheta^{\la}_{\mu,2}+
\sum_{r} rm_r(\mu)\,\vartheta^{\la}_{\mu_{\uparrow (r)}}=
 \al \sum_{i=1}^{l(\la)} \binom{\la}{\la_{(i)}} \,
\big(\la_i-1-(i-1)/\al\big)\,
\vartheta^{\la_{(i)}}_{\mu}.\quad
\end{equation}

\subsection{The operator $E_0$}

We have easily
\[E_0p_\mu=\sum_{r} rm_r(\mu)p_{\mu_{\downarrow (r)}}.\]
Specializing Proposition 4, and identifying coefficients of $N$ on both sides of the identity, we get
\begin{equation}
\vartheta^{\la}_{\mu}=
\sum_{i=1}^{l(\la)+1} c_i(\la) \, \vartheta^{\la^{(i)}}_{\mu},
\end{equation} 
which may be also obtained from Proposition 3 (ii). 

If we identify terms in $N^0$, we obtain
\begin{equation}
\sum^{\bullet}_{r} rm_r(\mu)\, \vartheta^{\la}_{\mu_{\downarrow (r)}}=
\sum_{i=1}^{l(\la)+1} c_i(\la) \, \big(\la_i-(i-1)/\al\big)\,
\vartheta^{\la^{(i)}}_{\mu}.
\end{equation}
Here for clarity of display, we use the symbol $\sum^{\bullet}$ to recall that for $r=2$ the partition $\mu_{\downarrow (2)}$ has one part $1$. Therefore due to (3.2), its contribution must be multiplied by the factor $n-|\mu_{\downarrow (2)}|+m_1(\mu_{\downarrow (2)})=n-k+2$. 

\subsection{The operator $D_1$}

We have $E_1p_\mu=kp_\mu$ by homogeneity and $(p_1p_\mu)^\#(\la)=(n-k)(p_\mu)^\#(\la)$ by Proposition 1 (i). 
Using (5.1) the identity of Proposition 6 then reads
\begin{multline*}
\sum_{r,s} rsm_r(\mu)(m_s(\mu)-\delta_{rs})\, 
\vartheta^{\la}_{\mu_{\downarrow (rs)}}
+(2N-1) \sum^{\bullet}_{r}rm_r(\mu)\, \vartheta^{\la}_{\mu_{\downarrow (r)}}\\
+\beta \sum^{\bullet}_{r} r(r-1)m_r(\mu) \, \vartheta^{\la}_{\mu_{\downarrow (r)}}
+\al \sum^{\bullet}_r rm_r(\mu)\sum_{i=1}^{r-2}\,
\vartheta^{\la}_{\mu_{\uparrow (i,r-i-1)}} 
+\big(n+k+N(N-1)/\al\big) \,
\vartheta^{\la}_{\mu}\\ = \al
\sum_{i=1}^{l(\la)+1} c_i(\la) \,\big(\la_i+(N-i)/\al\big)\,
\big(\la_i+(N-i+1)/\al\big)\,  \vartheta^{\la^{(i)}}_{\mu}.
\end{multline*}

As before, the symbol $\sum^{\bullet}$ is used to specify some particular cases :
\begin{enumerate} 
\item[(i)]The partition $\mu_{\downarrow (2)}$ and for $r\ge 4$ the partition $\mu_{\uparrow (1,r-2)}$ have one part $1$. Therefore due to (3.2), their contributions must be multiplied by $(n-k+2)$. 
\item[(ii)]The partition $\mu_{\uparrow (1,1)}$ has two parts $1$. Therefore its contribution must be multiplied by $(n-k+2)(n-k+3)/2$.
\end{enumerate}

Clearly if we identify the coefficients of $N^2$ on both sides, we recover (5.4). Then if we identify the coefficients of $N$, we recover (5.5). Finally the remaining terms give
\begin{multline}
\sum_{r,s} rsm_r(\mu)(m_s(\mu)-\delta_{rs})\, 
\vartheta^{\la}_{\mu_{\downarrow (rs)}}
\\ +\beta \sum^{\bullet}_{r} r(r-1)m_r(\mu) \, \vartheta^{\la}_{\mu_{\downarrow (r)}}
+\al \sum^{\bullet}_r rm_r(\mu)\sum_{i=1}^{r-2}\,
\vartheta^{\la}_{\mu_{\uparrow (i,r-i-1)}}\\
=-(n+k) \, \vartheta^{\la}_{\mu}+\al \,
\sum_{i=1}^{l(\la)+1} c_i(\la) \,
\big(\la_i-(i-1)/\al\big)^2\, \vartheta^{\la^{(i)}}_{\mu}.\qquad
\end{multline}

All coefficients $\vartheta^{\la}_{\rho}$ appearing in the left-hand side correspond to partitions $|\rho|=k-1$. Below this property will be crucial for our purpose.

\section{Rectangular shape}

We are now in a position to prove the following \textit{weak version} of our conjecture for $m=1$, i.e. when $\la=p \times q$, the rectangular shape formed by $p$ parts equal to $q$.
\begin{theor}
Let $\la=p \times q$ and $\mu$ a partition with $m_1(\mu)=0$ and $|\mu|=k\le |\la|=pq$. The quantity $(-1)^{k} \,z_{\mu}\, \theta^{\la}_{\mu,1^{pq-k}}$ is a polynomial in the indeterminates $(p,-q)$ and $\beta=\alpha-1$, with nonnegative rational coefficients.
\end{theor}

\begin{proof}
We shall use induction on the weight $|\mu|=k$. The property is verified for $k=2$. Actually it is well known (see \cite[p.384]{Ma}, ~\cite[p.106]{S0} or \cite[p.68]{La3}) that
\[2\theta^{\la}_{2,1^{n-2}}=2\al d_1(\la)= pq(\al q-p-\beta).\]

Obviously there are only two partitions $\la^{(i)}$, corresponding respectively to $i=1$ and $i=p+1$. We have
\begin{equation*}
\begin{split}
\la^{(1)}=(q+1,q,\ldots,q), \quad& \quad \la^{(p+1)}=(q,\ldots,q,1),\\
c_1(\la)=\frac{p}{p+\al q}, \quad& \quad c_{p+1}(\la)=\frac{\al q}{p+\al q}.
\end{split}
\end{equation*}
We consider the linear system formed by (5.4), (5.5) and (5.6). Firstly we evaluate the quantities $\vartheta^{\la^{(1)}}_{\mu}$ and $\vartheta^{\la^{(p+1)}}_{\mu}$ by solving (5.4) and (5.5). This easily yields 
\begin{equation}
\begin{split}
\vartheta^{\la^{(1)}}_{\mu}&=
\vartheta^{\la}_{\mu}+\frac{\al}{p}
\sum^{\bullet}_{r}rm_r(\mu)\, 
\vartheta^{\la}_{\mu_{\downarrow (r)}}
\\
\vartheta^{\la^{(p+1)}}_{\mu}&=
\vartheta^{\la}_{\mu}
-\frac{1}{q}\sum^{\bullet}_{r}rm_r(\mu)\, \vartheta^{\la}_{\mu_{\downarrow (r)}}.
\end{split}
\end{equation}
Secondly we insert both values in (5.6). The sum in the right-hand side writes
\[\al q^2 c_1(\la) \, \vartheta^{\la^{(i)}}_{\mu}
+p^2 c_{p+1}(\la) \, \vartheta^{\la^{(p+1)}}_{\mu}/\al
= pq \,\vartheta^{\la}_{\mu}
+(\al q-p) \sum^{\bullet}_{r}rm_r(\mu)\, 
\vartheta^{\la}_{\mu_{\downarrow (r)}},\]
so that finally (5.6) becomes
\begin{multline}
\sum_{r,s} rsm_r(\mu)(m_s(\mu)-\delta_{rs})\, 
\vartheta^{\la}_{\mu_{\downarrow (rs)}}
+\beta \sum^{\bullet}_{r} r(r-1) m_r(\mu) \, \vartheta^{\la}_{\mu_{\downarrow (r)}}\\
+(p-\al q ) \sum^{\bullet}_{r}rm_r(\mu)\, \vartheta^{\la}_{\mu_{\downarrow (r)}}
+\al \sum^{\bullet}_r rm_r(\mu)\sum_{i=1}^{r-2}\,
\vartheta^{\la}_{\mu_{\uparrow (i,r-i-1)}}
=-k \, \vartheta^{\la}_{\mu}.\quad
\end{multline}

This is an inductive formula expressing $k \,\vartheta^{\la}_{\mu}$ as an integral combination of $\vartheta^{\la}_{\rho}$ with $|\rho|=k-1$. Multiplying both sides by $(-1)^{k-1}$ and substituting $-q$ to $q$, we obtain that $(-1)^{k}  \vartheta^{\la}_{\mu}$ is a polynomial in the indeterminates $(p,-q,\beta)$, with \textit{nonnegative rational} coefficients.

Observe that the partitions having parts $1$ do not create any difficulty since their contribution is multiplied respectively by $-(-pq+k-2)$ for $\mu_{\downarrow (2)}$ or $\mu_{\uparrow (1,r-2)}$, and by $(-pq+k-2)(-pq+k-3)/2$ for $\mu_{\uparrow (1,1)}$. 
\end{proof}

We emphasize that the previous argument does not allow to conclude that the coefficients of $\vartheta^{\la}_{\mu}$ are \textit{integers}. It remains to prove that the coefficients of $k \, \vartheta^{\la}_{\mu}$ are \textit{divisible} by $k$. The recurrence (6.2) shows empirical evidence of this fact, but we are in lack of a proof. 

Below are the first steps of our recurrence, for $k\le 6$.
\begin{eqnarray*}
-\vartheta^{\la}_{2} & = & pq(p-\al q+\beta),\\
-\vartheta^{\la}_{3} & = & \vartheta^{\la}_{2}(p-\al q+2\beta)+\al pq(pq-1), \\
-\vartheta^{\la}_{4} & = & 
\vartheta^{\la}_{3}(p-\al q+3\beta)+2\al(pq-2)\vartheta^{\la}_{2},\\ 
-\vartheta^{\la}_{22} & = & 2\vartheta^{\la}_{3}+
(p-\al q+\beta)(pq-2)\vartheta^{\la}_{2},\\
-\vartheta^{\la}_{5} & = & \vartheta^{\la}_{4}(p-\al q+4\beta)+2\al(pq-3)\vartheta^{\la}_{3}+\al \vartheta^{\la}_{22},\\
-5\vartheta^{\la}_{32} & = & 12\vartheta^{\la}_{4}+
2(p-\al q+\beta)(pq-3)\vartheta^{\la}_{3} \\& & \
+3(p-\al q+2\beta)\vartheta^{\la}_{22}
+3\al(pq-2)(pq-3)\vartheta^{\la}_{2},\\
-\vartheta^{\la}_{6} & = & \vartheta^{\la}_{5}(p-\al q+5\beta)+2\al(pq-4)\vartheta^{\la}_{4}
+2\al \vartheta^{\la}_{32},\\
-6\vartheta^{\la}_{42} & = & 16\vartheta^{\la}_{5}+
2(p-\al q+\beta)(pq-4)\vartheta^{\la}_{4} \\& & \
+4(p-\al q+3\beta)\vartheta^{\la}_{32}
+8\al(pq-4)\vartheta^{\la}_{22},\\
-\vartheta^{\la}_{33} & = & 3\vartheta^{\la}_{5}+
(p-\al q+2\beta)\vartheta^{\la}_{32} 
+\al(pq-3)(pq-4)\vartheta^{\la}_{3},\\
-\vartheta^{\la}_{222} & = & 4\vartheta^{\la}_{32}+
(p-\al q+\beta)(pq-4)\vartheta^{\la}_{22}.
\end{eqnarray*}

\begin{prop}
Let $\Lambda$ be a partition obtained by adding or substracting one node to the rectangular shape $\la=p \times q$. The quantity $(-1)^{k} \, \vartheta^{\Lambda}_{\mu}$ is a polynomial in the indeterminates $(p,-q,\beta)$, with nonnegative rational coefficients.
\end{prop}
\begin{proof}
From the recurrence formula (6.2), it is clear that any $ \vartheta^{\la}_{\rho}$ is divisible by $pq$. Thus by (6.1), the assertion is true for both partitions $\la^{(1)}=(q+1,q,\ldots,q)$ and $\la^{(p+1)}=(q,\ldots,q,1)$.
Obviously there is only one partition $\la_{(i)}$, corresponding to $i=p$. We have $\la_{(p)}=(q,\ldots,q,q-1)$ and
$\binom{\la}{\la_{(p)}}=pq$.
Therefore (5.2) writes
\begin{equation}
(pq-k) \vartheta^{\la}_{\mu} = pq \, \vartheta^{\la_{(p)}}_{\mu}.
\end{equation}
Hence the statement for $\la_{(p)}$.
\end{proof}

\section{Final remark} 

It is worth showing why the proofs given for $\alpha=1$~\cite{Ra,S1} are difficult to extend when $\al$ is arbitrary. This will also produce a non trivial result.

The proof given in~\cite{Ra} starts from the following formula, proved in~\cite[(15.21)]{Ok2}:
\[\widehat{\chi}^\la_{\mu,1^{n-k}}=\sum_{|\rho|=k} s_\rho^\dag(\la)
\, \chi^\rho_\mu.\]
Here $s_\rho^\dag$ is an appropriate normalization of the shifted Schur function, i.e. the shifted Jack polynomial corresponding to $\al=1$. Hence $s_\rho^\dag(\la)$ is, up to some normalization, the generalized binomial coefficient $\binom{\la}{\rho}$ for $\al=1$. 

Due to (1.1) this can be rewritten under the form
\[\binom{n-k+m_1(\mu)}{m_1(\mu)} \theta^{\la}_{\mu,1^{n-k}}=
\sum_{|\rho|=k}\binom{\la}{\rho} \theta^{\rho}_{\mu}.\]
It is a remarkable fact that this property keeps true when $\al$ is arbitrary. Actually it is exactly (3.3) written for $r=k$.

Moreover for $\la=p \times q$ the binomial coefficients $\binom{\la}{\rho}$ are \textit{explicitly known} for any $\rho$. They have been computed, in the more general context of Macdonald polynomials, in~\cite[Theorem 11, p.313]{La2}. As a limit case~\cite[p.321]{La2} we have
\[\binom{p \times q}{\rho}=(-\al)^{|\rho|}(-q)_\rho \, J_\rho^\sharp(1^p)=
(-1)^{|\rho|}\al^{2|\rho|}\,(-q)_\rho\,\frac{(p/\al)_\rho}{j_\rho},\]
the last equality being a consequence of (2.1).

Thus for $\la=p \times q$ we have
\[\binom{pq-k+m_1(\mu)}{m_1(\mu)} \theta^{p \times q}_{\mu,1^{pq-k}}=(-1)^k\al^{2k}
\sum_{|\rho|=k} (p/\al)_\rho\,(-q)_\rho\,\frac{\theta^{\rho}_{\mu}}{j_\rho}.\]
And Theorem 1 appears equivalent to the following result, which seems difficult to prove directly.
\begin{theor}
Let $p,q$ be two indeterminates. For any partition $\mu$ the quantity
\[\al^{2|\mu|-1}\, z_{\mu}\,\sum_{|\rho|=|\mu|}(p)_\rho(q)_\rho\frac{\theta^{\rho}_{\mu}}{j_\rho} \]
is a polynomial in $p$, $q$, $\beta$ with nonnegative rational coefficients.
\end{theor}

We conjecture (i) the integrality of these coefficients, (ii) the existence of some natural expression as a polynomial in $(\al,\beta)$. However we have no conjectured combinatorial interpretation for such an expression. 

As an example for $\mu=(3,2)$ we have 
\begin{multline*}
6\al^{9}\,\sum_{|\rho|=5}(p)_\rho(q)_\rho\frac{\theta^{\rho}_{32}}{j_\rho}=
pq \Big (
(p^4q +4p^3q^2 +4p^2q^3+pq^4 ) \al^4
 +(4p^3q +9p^2q^2+4pq^3) \al^3\beta\\ +5(p^2q+pq^2) \al^2\beta^2  
+2pq \al\beta^3 +(6p^3+31p^2q+31pq^2+6q^3)\al^3 \\ 
+(30p^2+79pq+30q^2) \al^2\beta +48(p+q)\al\beta^2 +24\beta^3
+18(p+q)\al^2 +24 \al\beta\Big).
\end{multline*} 

\section{The role of $\beta$}

Our results for $\la=p \times q$ suggest the existence of some unknown underlying pattern where $\beta=\al-1$ would play a role as important as $\al$. Actually the recurrence formula (6.2) gives an expression of the quantities $\vartheta^{\la}_{\rho}$ as ``positive'' polynomials in both parameters $\al$ and $\beta$.

As mentioned in the introduction, we conjecture that this property is general, i.e. that for any $\la=\bmp \times \bmq$ the quantities $\vartheta^{\la}_{\rho}$ have some natural expression as polynomials in $(\alpha,\beta)$ with nonnegative integer coefficients.

Actually we can say more. In this paper $\beta=\al-1$ was of course never considered as being independent of $\al$. However for $\la=p \times q$, there is a strong empirical evidence that $\al$ and $\beta$ might be considered as two \textit{independent} parameters. 

This strange fact can be seen on the identities obtained by our method. These identities, satisfied for any $\la$, are expressed in terms of $\al$ and $\beta=\al-1$. But at least for $\la=p \times q$, this restriction may be dropped. Here is a typical example, among many others.

\subsection{Example}

Let $\la$ be arbitrary and $\mu$ a partition with $m_1(\mu)=0$ and $|\mu|=k\le |\la|=n$. Writing Proposition 8 for $j=2$, we obtain
\begin{multline*}
D_2p_\mu=p_\mu \Big(\frac{2}{\al}(N-1)|\mu|+
\frac{\beta}{\al} \sum_{r} r(r-1)m_r(\mu)\\
+\sum_{r, s} rsm_r(\mu)(m_s(\mu)-\delta_{rs}) \frac{p_{r+s}}{p_rp_s}+\frac{1}{\al}\sum_r rm_r(\mu) \sum_{i=1}^{r-1} \frac{p_ip_{r-i}}{p_r}\Big).
\end{multline*}

For any integers $r,s\ge 1$ we denote $\mu_{\Downarrow (rs)}$ the partition (if it exists) obtained from $\mu$ by substracting one part $r$ and one part $s$, and adding one part $r+s$. Similarly we write $\mu_{\Uparrow (rs)}$ the partition (if it exists) obtained from $\mu$ by adding one part $r$ and one part $s$, and substracting one part $r+s$. We have $|\mu_{\Downarrow (rs)}|=|\mu_{\Uparrow (rs)}|=k$, $l(\mu_{\Downarrow (rs)})=l(\mu)-1$ and $l(\mu_{\Uparrow (rs)})=l(\mu)+1$. 

Then we obtain
\begin{multline*}
\Big(D_2^\dag -\frac{\beta}{\al} \sum_{r} r(r-1)m_r(\mu)\Big)p_\mu
=\\
\sum_{r, s} rsm_r(\mu)(m_s(\mu)-\delta_{rs}) p_{\mu_{\Downarrow (rs)}}+
\frac{1}{\al}\sum_r rm_r(\mu)\sum_{i=1}^{r-1}p_{\mu_{\Uparrow (i,r-i)}},
\end{multline*}
and Proposition 7, specialized for $f=p_\mu$, writes
\begin{multline*}
\sum_{r,s} rsm_r(\mu)(m_s(\mu)-\delta_{rs})\, 
\vartheta^{\la}_{\mu_{\Downarrow (rs)}}
+\al \sum^{\bullet} rm_r(\mu)\sum_{i=1}^{r-1}\,
\vartheta^{\la}_{\mu_{\Uparrow (i,r-i)}}\\
+2\sum_{r} rm_r(\mu)\,\vartheta^{\la}_{\mu_{\uparrow (r)}}
+ \vartheta^{\la}_{\mu,2}
= \Big(2\al d_1(\la)-\beta \sum_{r} r(r-1)m_r(\mu)\Big) \vartheta^{\la}_{\mu}.
\end{multline*}
Here, as before, the symbol $\sum^{\bullet}$ is used to specify some particular cases :
\begin{enumerate} 
\item[(i)]For $r\ge 3$ the partition $\mu_{\Uparrow (1,r-1)}$ has one part $1$. Therefore due to (3.2), its contribution must be multiplied by $(n-k+1)$.
\item[(ii)]The partition $\mu_{\Uparrow (1,1)}$ has two parts $1$ and its contribution must be multiplied by a factor $(n-k+1)(n-k+2)/2$.
\end{enumerate} 
 
This identity is true for any $\la$ and for any $\mu$ with $m_1(\mu)=0$. But if we write it for $\la=p \times q$ and replace each $\vartheta^{\la}_{\rho}$ by its value inductively defined through (6.2), we obtain an identity involving \textit{two independent parameters} $\al$ and $\beta$.

The simplest case is given by $\mu=(2)$, i.e.
\[2(\al d_1(\la)-\beta) \vartheta^{\la}_{2}=2\al n(n-1)+4\vartheta^{\la}_{3}+\vartheta^{\la}_{22}.\]
Once written for $\la=p \times q$, this becomes 
\begin{multline*}
(pq(\al q-p-\beta))^2 =\\ 2\beta pq(\al q-p-\beta) 
+2pq(\al q-p-\beta)(\al q-p-2\beta)
+pq(pq-2)(\al q-p-\beta)^2,
\end{multline*}
which is an identity in four independent variables $(p,q,\al,\beta)$.
The reader may check the fact for other values of $\mu$, for instance for $\mu=(3,2)$ where the identity reads 
\begin{equation*}
(2\al d_1(\la)-8\beta)\,\vartheta^{\la}_{32} = 
12\vartheta^{\la}_{5}+2\al (n-3)(n-4)\vartheta^{\la}_{3}
+6\al(n-4)\vartheta^{\la}_{22} 
+6\vartheta^{\la}_{42}
+4\vartheta^{\la}_{33}+\vartheta^{\la}_{322}.
\end{equation*}

We emphasize that this example is far from being isolated. Actually it seems that for $\la=p \times q$ we may consider $(\al,\beta)$ as independent parameters in any identity linking the quantities $\vartheta_\rho^\la$.

\subsection{Extension of Jack polynomials}

It is thus a reasonable attempt to drop the condition $\beta=\al-1$, and to generalize Jack polynomials as follows. Let $\al,\beta$ be two independent parameters. For each rectangular shape $\la=p \times q$, we set
\begin{equation*}
J_{\la}(\al,\beta)=\sum_{\mu,\, m_1(\mu)=0} \vartheta^{\la}_{\mu}(\al,\beta) \,\frac{p_{\mu}}{z_{\mu}}p_1^{|\la|-|\mu|},
\end{equation*}
with $\vartheta^{\la}_{\mu}(\al,\beta)$ inductively defined by (6.2). 

Because of Proposition 9, this extension can also be defined when $\la$ is obtained by adding or substracting one node to a rectangular shape $p \times q$. Namely,
\begin{enumerate} 
\item[(i)]For $\la=(q+1,q,\ldots,q)$ and $\la=(q,\ldots,q,1)$, we define $\vartheta^{\la}_{\mu}(\al,\beta)$ by (6.1).
\item[(ii)]For $\la=(q,\ldots,q,q-1)$ we define $\vartheta^{\la}_{\mu}(\al,\beta)$ by (6.3).
\end{enumerate} 

However something new here happens. For $\beta=\al-1$ the quantities $\vartheta^{\la}_{\mu}(\al,\beta)$ inductively defined by (6.2) vanish for $|\mu|>|\la|$, which reflects the obvious fact that the sum must be restricted to $|\mu|\le |\la|$.

Remarkably this vanishing property is no longer satisfied for arbitrary $\beta$. By induction it may be seen that for $|\mu|>|\la|$, any $\vartheta^{\la}_{\mu}(\al,\beta)$ is divisible by $(\al-\beta-1)$. Here are a few examples. For $\la=(1)$, i.e. $p=q=1$, we have
\[\vartheta^{(1)}_{2}=  (\al-\beta-1), \quad \quad
\vartheta^{(1)}_{3}=  
(\al-\beta-1)(\al-2\beta-1).\]
And for $\la=1^3$, i.e. $p=3, q=1$, we have
\[\vartheta^{1^3}_{32} =  
-9(\al-\beta-1)(\al-2\beta-3)(\al-3\beta-9).\]
Does it mean that $J_{\la}(\al,\beta)$ is actually \textit{a series}, terminating for $\beta=\al-1$ ?

Many other questions are still open. Have these objects a more natural definition ? Are they eigenvectors of some differential operators ? Are they orthogonal for some scalar product ? 

We conjecture that such a two-parameters extension $J_\la(\al,\beta)$ may be defined for any $\la= \bmp \times \bmq$.


\begin{thebibliography}{29}

\bibitem{F} V.\ F\'eray, \emph{Proof of Stanley's conjecture about irreducible character values of the symmetric group}, math.CO/0612090.
\bibitem{Ka}
J.\ Kaneko, \emph{Selberg integrals and hypergeometric functions associated with Jack polynomials}, SIAM J. Math. Anal. \textbf{24} (1993), 1086--1110.
\bibitem{KS}
F.\ Knop, S.\ Sahi, \emph{Difference equations and symmetric polynomials defined by their zeros}, Int. Math. Res. Not. \textbf{1996} (10), 473--486.
\bibitem{KS2}
F.\ Knop, S.\ Sahi, \emph{A recursion and a combinatorial formula for Jack polynomials}, Invent. Math. \textbf{128} (1997), 9--22.
\bibitem{La0}
M.\ Lassalle, \emph{Une formule de Pieri pour les polyn\^omes de Jack}, C. R. Acad. Sci. Paris, S\'er. I Math., \textbf{309} (1989), 941--944.
\bibitem{La1}
M.\ Lassalle, \emph{Une formule du bin\^ome g\'en\'eralis\'ee pour les 
polyn\^omes de Jack}, C. R. Acad. Sci. Paris, S\'er. I Math., \textbf{310} (1990).
\bibitem{La2}
M.\ Lassalle, \emph{Coefficients binomiaux g\'en\'eralis\'es et 
polyn\^omes de Macdonald}, J. Funct. Anal. \textbf{158} (1998), 
289--324.
\bibitem{La3}
M.\ Lassalle, \emph{Some combinatorial conjectures for Jack polynomials}, Ann. Comb. \textbf{2} (1998), 61--83.
\bibitem{Ma}
I.\ G.\ Macdonald, \emph{Symmetric functions and Hall polynomials}, Clarendon Press, second edition, Oxford, 1995.
\bibitem{Ok1}
A.\ Okounkov, G.\ Olshanski, \emph{Shifted Jack polynomials, binomial formula and applications}, Math. Res. Lett. \textbf{4} (1997), 69--78.
\bibitem{Ok2}
A.\ Okounkov, G.\ Olshanski, \emph{Shifted Schur functions}, St. Petersburg Math. J. \textbf{9} (1998), 239--300.
\bibitem{Ok3}
A. Okounkov, \emph{(Shifted) Macdonald polynomials, $q$-integral representation and combinatorial formula}, Compositio. Math. \textbf{112} (1998), 147--182.
\bibitem{Ra} A.\ Rattan, \emph{Positivity results for Stanley's character polynomials}, Journal of Algebra, to appear, math.RT/0601186.
\bibitem{S0}
R.\ P.\ Stanley, \emph{Some combinatorial properties of Jack
symmetric functions}, Adv.\ Math., \textbf{77} (1989), 76--115.
\bibitem{S1} R.\ P.\ Stanley, \emph{Irreducible symmetric group characters of rectangular shape}, S\'em. Lothar. Combin., \textbf{50} (2003), article B50d.
\bibitem{S2} R.\ P.\ Stanley, \emph{A conjectured combinatorial interpretation of the normalized irreducible character values of the symmetric group}, math.CO/0606467.
\bibitem{W}\texttt{http://igm.univ-mlv.fr/{\textasciitilde}lassalle/conj.html}


\end{thebibliography}
\end{document}